\documentclass[11pt]{amsart}
\usepackage{amsmath,amssymb,vmargin,amsthm}
\usepackage[latin1]{inputenc}
\setpapersize{A4}
\setmarginsrb{26mm}{2cm}{26mm}{2cm}{0cm}{1cm}{0cm}{1cm}

\newcommand\N{{\mathbb N}}

\newcommand\R{{\mathbb R}}

\newcommand\T{{\mathcal T}}
\newcommand\Z{{\mathbb Z}}
\newcommand\set[2]{\left\{#1\left|\,#2\right.\right\}}

\newcommand\aut{\operatorname{Aut}}

\newcommand\abs[1]{|#1|}

\newcommand\quo[2]{#1\bmod#2}
\newcommand\gen[1]{\left\langle#1\right\rangle}
\newcommand\stab[2][]{{\operatorname{Stab}}_{#1}(#2)}

\newcommand\st[2]{\quo{#1}{#2}}
\newcommand\sym[1]{\operatorname{Sym}(#1)}

\newtheorem{theorem}{Theorem}[section]
\newtheorem{proposition}[theorem]{Proposition}
\newtheorem{corollary}[theorem]{Corollary}
\newtheorem{lemma}[theorem]{Lemma}
\theoremstyle{remark}
\newtheorem{remark}[theorem]{Remark}
\newtheorem{example}[theorem]{Example}
\theoremstyle{definition}
\newtheorem{definition}[theorem]{Definition}

\newcommand\diminf{\underline{\dim}}
\def\pair<#1>{{\ll}#1{\gg}}

\title{Hausdorff dimension of some groups acting on the binary tree}
\author{Olivier Siegenthaler}
\address{EPFL SB IMB MAD, Station 8, CH-1015 Lausanne, Switzerland}
\email{olivier.siegenthaler@epfl.ch}
\date{October 1, 2007}

\begin{document}
\bibliographystyle{habbrv}

\begin{abstract}
Based on the work of Abercrombie~\cite{abercrombie:sspr}, Barnea and Shalev~\cite{barnea_shalev:hd} gave an explicit formula for the Hausdorff dimension of a group acting on a rooted tree. We focus here on the binary tree $\T$. Ab\'ert and Vir\'ag~\cite{abert_virag:drgart} showed that there exist finitely generated (but not necessarily level-transitive) subgroups of $\aut\T$ of arbitrary dimension in $[0,1]$.

In this article we explicitly compute the Hausdorff dimension of the level-transitive spinal groups. We then show examples of $3$-generated spinal groups which have transcendental Hausdroff dimension, and exhibit a construction of $2$-generated groups whose Hausdorff dimension is $1$.
\end{abstract}
\maketitle

\footnotetext{The author gratefully acknowledges support from the Swiss National Fund for Scientific Research (grant 105469/1).}

\section{Introduction}
Although it is known~\cite{abert_virag:drgart} that finitely generated subgroups of $\aut\T$ may have arbitrary Hausdorff dimension, there are only very few explicit computations in the literature. Further, all known examples have rational dimension, starting with the ``first Grigorchuk group" whose dimension is $5/8$, as we will see below. In this article we give two explicit constructions of finitely generated groups. We obtain groups of dimension $1$ on one hand, and groups whose dimension is transcendental on the other hand. This is achieved by computing the dimension of the so-called spinal groups acting on the binary tree, which are generalizations of the Grigorchuk groups.

We begin by recalling the definition of Hausdorff dimension in the case of groups acting on the binary tree. In Section~\ref{sec_sg}, we define the spinal groups we are interested in. Section~\ref{sec_main} is devoted to the statement of Theorem~\ref{th_main} which gives a formula for computing the Hausdorff dimension of any spinal group acting level-transitively on the binary tree; the proof is deferred to Section~\ref{sec_pr}. Finally, the construction of $3$-generated spinal groups of irrational Hausdorff dimension is given in Section~\ref{sec_irr}, and groups with Hausdorff dimension $1$ are exhibited in Section~\ref{sec_ful}.

\subsection*{Acknowledgements}
I wish to thank Laurent Bartholdi who suggested these questions to me as well as many ideas for answering them.

\section{Hausdorff Dimension}\label{sec_hd}
Let $\T$ be the infinite binary rooted tree, and let $\aut\T$ denote its automorphism group. It is known~\cite{barnea_shalev:hd} that the Hausdorff dimension of a closed subgroup $G$ of $\aut\T$ is
\begin{align*}
  \dim_H G=\liminf_{m\to\infty}\frac{\log\abs{\quo{G}{m}}}{\log\abs{\st{\aut \T}{m}}},
\end{align*}
where the ``mod $m$" notation stands for the action on the first $m$ levels of the tree (i.e. $(\quo{G}{m})=G/\stab[G]{m}$, where $\stab[G]{m}$ is the fixator of the $m$th level of the tree). Moreover, one easily computes $\abs{\st{\aut \T}{m}}=2^{2^m-1}$. This yields the more explicit formula
\begin{align*}
  \dim_H G=\liminf_{m\to\infty}\frac{\log_2\abs{\quo{G}{m}}}{2^m}.
\end{align*}
Below we will identify the vertices of $\T$ with the set of finite words over the alphabet $X=\{0,1\}$. We recall that there is a canonical decomposition of the elements $g\in\aut\T$ as
\begin{align*}
  g=\pair<g@0,g@1>\sigma,
\end{align*}
with $g@x\in\aut\T$ and $\sigma\in\sym X$. We will often identify $\sym X$ with $C_2$, the cyclic group of order 2, or with the additive group of $\Z/2\Z$, the finite field with 2 elements.

\section{Spinal groups}\label{sec_sg}
We only deal here with a specific case of the more general definition of spinal groups which can be found in~\cite{bartholdi_sunic:wpg}.

Before defining the spinal group $G_\omega$, we need a \emph{root group} $A$ which we will always take to be $A=\gen{a}=C_2$, and a \emph{level group} $B$, which will be the $n$-fold direct power of $C_2$. We think of $B$ as an $n$-dimensional vector space over $\Z/2\Z$.
Let $\omega=\omega_1\omega_2\dots$ be a fixed infinite sequence of non-trivial elements from $B^*$, the dual space of $B$, and let $\Omega$ denote the set of such sequences.
We let the non-trivial element $a$ of $A$ act on $\T$ by exchanging the two maximal subtrees. Next, we let each element $b\in B$ act via the recursive formul\ae\  $b=\pair<\omega_1(b),b@1>$ and $b@1^k=\pair<\omega_{k+1}(b),b@1^{k+1}>$. The \emph{spinal group} $G_\omega$ is the group generated by $A\cup B$.
Note that we required each $\omega_i$ to be non-trivial, and this implies that $G_\omega$ is level-transitive (i.e. $G_\omega$ acts transitively on $X^n$ for all $n\in\N$).

The \emph{syllable form} of $\omega$ is $\omega_1^{a_1}\omega_2^{a_2}\dots$ where $\omega_i\neq\omega_{i+1}$ and the $a_i$ denote multiplicities. In contrast to this we say that $\omega=\omega_1\omega_2\dots$ is in \emph{developed form}. We designate by $s_k$ the sum of the $k$ first terms of the sequence $(a_i)$, i.e. the length of the $k$-syllable prefix of $\omega$.

\section{Main theorem}\label{sec_main}

Before stating the main theorem, we need a few more definitions.
\begin{definition}
The \emph{shift} $\sigma:\Omega\to\Omega$ is defined on a sequence $\omega=\omega_1\omega_2\dots$ by $\sigma\omega=\omega_2\omega_3\dots$.
\end{definition}
\begin{definition}
For a positive integer $m$ and a sequence $\omega=\omega_1^{a_1}\omega_2^{a_2}\dots$ in syllable form, we define $s^{-1}(m)$ as the only integer that satisfies $s_{s^{-1}(m)-1}+1<m\leq s_{s^{-1}(m)}+1$.
\end{definition}
In other words, $s^{-1}(m)$ is the number of syllables of the prefix of $\omega$ of length $m-1$. If we write $\omega=\omega_1^{a_1}\omega_2^{a_2}\dots$ in syllable form, then $\{\omega_1,\dots,\omega_{s^{-1}(m)}\}$ is the set of elements of $B^*$ which appear in $(\quo{G}{m})$.
\begin{definition}
Let $\omega=\omega_1^{a_1}\omega_2^{a_2}\dots$ be in syllable form. Given $m\in\N$, we define $\dim(\quo{\omega}{m})$ as the dimension of the vector space spanned by $\omega_1,\dots,\omega_{s^{-1}(m)}$. We also define $\diminf(\omega)$ by
\begin{align*}
\diminf(\omega)=\liminf_{n\to\infty}\left\{\lim_{m\to\infty}\dim(\quo{\sigma^n\omega}{m})\right\}.
\end{align*}
\end{definition}

The next theorem expresses the Hausdorff dimension of any level-transitive spinal group acting on the binary tree. Its proof is given in Section~\ref{sec_pr}.
\begin{theorem}\label{th_main}
Consider $\omega=\omega_1^{a_1}\omega_2^{a_2}\dots$ in syllable form with $\diminf(\omega)=n\geq2$. The Hausdorff dimension of $G_\omega$ is given by
\begin{align*}
  \dim_HG_{\omega}=\frac{1}{2}\liminf_{k\to\infty}\left(\frac{\Sigma_k}{2^{s_k}}
+\frac{1}{2^{s_k}}\sum_{i=2}^{n-1}2^{s_{\lambda_i}}\left(2-\frac{1}{2^{a_{\lambda_i}}}\right)
+\frac{1}{2^{a_k}}\left(1-\frac{1}{2^{a_{k-1}}}\right)\right),
\end{align*}
where $\Sigma_k=\sum_{j=1}^{k}2^{s_{j-1}}a_j$, and for each $i\in\{2,\dots,n-1\}$ we let $\lambda_i(k)$ be the smallest integer such that
\begin{align*}
\dim(\quo{\sigma^{s_{\lambda_i(k)}}\omega}{(s_{k+1}-s_{\lambda_i(k)})})=i.
\end{align*}
\end{theorem}

\begin{remark}
In case the sequence $\omega$ is eventually constant, i.e. if $\diminf(\omega)=1$, then Proposition~\ref{prop:sizegmodm} can be used to show that $\log_2\abs{\quo{G_\omega}{m}}$ grows linearly with $m$, whence $\dim_HG_\omega=0$.
\end{remark}

In the special case where $B=C_2\times C_2$ and $\omega$ is not eventually constant, one can use the following corollary to compute $\dim_HG_\omega$:

\begin{corollary}\label{cor_dim2}
If $\diminf(\omega)=2$, then
\begin{align*}
  \dim_HG_{\omega}=\frac{1}{2}\liminf_{k\to\infty}\left(\frac{\Sigma_k}{2^{s_k}}
+\frac{1}{2^{a_k}}\left(1-\frac{1}{2^{a_{k-1}}}\right)\right).
\end{align*}
\end{corollary}

\begin{example}
Consider $B=C_2\times C_2$. Let $\omega_1,\omega_2,\omega_3$ be the three non-trivial elements of $B^*$, and consider $\omega=\omega_1\omega_2\omega_3\omega_1\omega_2\omega_3\dots$. The group $G_\omega$ is the ``first Grigorchuk group", first introduced in~\cite{grigorchuk} (see also~\cite{aleshin}). Since the corresponding integer sequence is given by $a_k=1$ for all $k$, then $s_k=k$ and $\Sigma_k=2^k-1$. The last corollary yields
\begin{align*}
  \dim_HG_{\omega}=\frac{5}{8}.
\end{align*}
\end{example}

\begin{example}
The Hausdorff dimension of some spinal groups is computed in~\cite{sunic:hdfsg}. Consider $B=C_2^n$, and fix a functional $\phi\in B^*$ and an automorphism $\rho$ (i.e. an invertible linear transformation) of $B$. We consider the sequence $\omega=\omega_1\omega_2\dots$ defined by $\omega_1=\phi$ and $\omega_n=\rho^*(\omega_{n-1})$ for $n>1$, where $\rho^*$ denotes the adjoint automorphism of $\rho$. We restrict attention to the case where $\rho$ and $\phi$ are such that $\diminf(\omega)=n$ (this is a rephrasing of the condition in~\cite{sunic:hdfsg}, Proposition~2). This implies that every sequence of $n$ consecutive terms of $\omega$ generates $B^*$. Each syllable of $\omega$ has length $1$, so we can apply Theorem~\ref{th_main} with $a_i=1$ and $s_i=i$ for all $i$, and $\lambda_i(k)=k-i$ for $i\in\{2,\dots,n-1\}$ and $k\geq n$. This yields
\begin{align*}
  \dim_HG_{\omega}=1-\frac{3}{2^{n+1}},
\end{align*}
as was found in~\cite{sunic:hdfsg}.
\end{example}

\section{Construction of finitely generated groups of irrational Hausdorff dimension}\label{sec_irr}

Throughout this section we restrict to the case $\diminf(\omega)=2$.
Let $D\subset[0,1]$ be the set of possible Hausdorff dimension for $G_\omega$.
More precisely, $D$ is defined as
\begin{align*}
D=\set{\lambda\in\R}{\exists\omega\in\Omega\text{ with }\diminf(\omega)=2\text{ and }\dim_HG_\omega=\lambda}.
\end{align*}
Although it is not an easy thing to tell whether a given $x\in[0,1]$ lies in $D$, we are able to show the following. Let $C$ denote the Cantor set constructed by removing the second quarter of the unit interval, and iterating this process on the obtained intervals. It is easy to see that $C$ is compact and totally disconnected, and contains transcendental elements.

\begin{theorem}\label{th_CinD}
The set $D$ contains several copies of $C$, each one being the image of $C$ under an affine map with rational coefficients.
\end{theorem}
\begin{corollary}
The set $D$ contains transcendental elements.
\end{corollary}

\begin{proof}[Proof of Theorem~\ref{th_CinD}]
We first define the functions
\begin{align*}
f_{a,n}:\R&\to\R\\
x&\mapsto\frac{x}{2^a}+\frac{a+n}{2^a}.
\end{align*}
We will only pay attention to the case where $a$ is a strictly positive integer and $n\in\Z$, and we will simply write $f_a$ instead of $f_{a,0}$.

Consider $\omega=\omega_1^{a_1}\omega_2^{a_2}\dots$. Starting from Corollary~\ref{cor_dim2}, we can see that
\begin{align}
\dim_HG_\omega&=\frac{1}{2}\liminf_{k\to\infty}\left(\frac{\Sigma_k}{2^{s_k}}+\frac{1}{2^{a_k}}\left(1-\frac{1}{2^{a_{k-1}}}\right)\right)\nonumber\\
&=\frac{1}{2}\liminf_{k\to\infty}(f_{a_k,1}\circ f_{a_{k-1},-1}\circ f_{a_{k-2}}\circ\cdots\circ f_{a_1})(0).\label{eq_sa}
\end{align}
Observe that the functions $f_1$ and $f_2$ define an iterated function system whose invariant set $C'$ is the image of $C$ under the map $x\mapsto\frac{x+2}{3}$. Indeed $x_1=1$ and $x_2=\frac{2}{3}$ are the fixed points of $f_1$ and $f_2$ respectively. Write $\Delta=x_1-x_2$. Then $f_1([x_2,x_2+\Delta])=[x_2+\frac{\Delta}{2},x_2+\Delta]$ and $f_2([x_2,x_2+\Delta])=[x_2,x_2+\frac{\Delta}{4}]$.

Now fix a point $\widehat{x}\in C'$. There exists a sequence $(b_i)\in\{1,2\}^\N$ such that
\begin{align}
\widehat{x}=\lim_{k\to\infty}(f_{b_0}\circ\cdots\circ f_{b_k})(y)\label{eq_sb}
\end{align}
for any point $y\in\R$. We call the sequence $(b_0,b_1,\dots)$ the \emph{code} of $\widehat{x}$.
Notice that the main differences between~\eqref{eq_sa} and~\eqref{eq_sb} are the ordering of the factors, and the limit which is a $\liminf$ in~\eqref{eq_sa}.

Fix $s\in\N$, $s>2$. We define the sequence $(a_i)\in\{1,2,s\}^\N$ by
\begin{align*}
(a_1,a_2,\dots) = (b_0,s,b_1,b_0,s,b_2,b_1,b_0,s,\dots)
\end{align*}
The sequence $(a_i)$ thus consists of prefixes of the sequence $(b_i)$ of increasing length, written backwards, and separated by $s$.
We set $\omega=\omega_1^{a_1}\omega_2^{a_2}\omega_1^{a_3}\omega_2^{a_4}\dots$.
We will show that $\dim_HG_\omega$ is the image of $\widehat{x}$ under an affine map with rational coefficients.
Recall that $\dim_HG_\omega$ is given by the $\liminf$ of
\begin{align}\label{eq:liminf}
\frac{1}{2}(f_{a_k,1}\circ f_{a_{k-1},-1}\circ f_{a_{k-2}}\circ\cdots\circ f_{a_1})(0).
\end{align}
Notice that the maps $f_{a,n}$ are order-preserving, and observe that $(f_{\alpha,1}\circ f_{\beta,-1})(0)\geq\frac{3}{4}$ if $\alpha\in\{1,2\}$ and $b\geq1$, while $(f_{\alpha,1}\circ f_{\beta,-1})(1)<\frac{3}{4}$ if $\alpha=s$ and $\beta\geq1$.
Therefore the lowest values in~\eqref{eq:liminf} are attained when $a_k=s$.
We conclude that
\begin{align*}
\dim_HG_\omega=\frac{1}{2}\lim_{k\to\infty}(f_{s,1}\circ f_{b_{0,-1}}\circ f_{b_1}\circ\cdots\circ f_{b_k})(0)=\frac{1}{2}(f_{s,1}\circ f_{b_0,-1}\circ f_{b_0}^{-1})(\widehat{x}).
\end{align*}

It should be noted that the maps $f_{s,1}\circ f_{b_0,-1}\circ f_{b_0}^{-1}$ do not have disjoint images. Nevertheless, it can be checked that
\begin{align*}
f_{1,-1}\circ f_{1}=f_{2,-1}.
\end{align*}
This implies that a point $x\in C'$ whose code is $2b_1b_2\dots$ is mapped to the same point as the point whose code is $11b_1b_2\dots$. On the other hand, the set of points whose code starts with a $1$ is just $f_1(C')$. Under the maps $f_{s,1}\circ f_{b_0,-1}\circ f_{b_0}^{-1}$, the set $f_1(C')$ is sent to $(f_{s,1}\circ f_{1,-1})(C')$, and the maps $f_{s,1}\circ f_{1,-1}$ are affine with rational coefficients and have disjoint images for all $s>2$. Thus $D$ contains a countable infinity of disjoint copies of $C'$.
\end{proof}

\section{Construction of full-dimensional finitely generated groups}\label{sec_ful}

We begin with a few easy statements which will be useful.
\begin{proposition}\label{prop_dima}
Let $H\leq G$ be subgroups of $\aut X^*$. Then $\dim_HH\leq\dim_HG$. Moreover, $\dim_HH=\dim_HG$ if the index of $H$ in $G$ is finite.
\end{proposition}
\begin{proof}
$H\leq G$ implies $\abs{\quo{H}{n}}\leq\abs{\quo{G}{n}}$ for all $n\geq0$. Thus $\dim_HH\leq\dim_HG$. Moreover, if $k=\abs{G:H}$ is finite then $k\abs{\quo{H}{n}}\geq\abs{\quo{G}{n}}$ for all $n\geq0$. This yields the second claim.
\end{proof}
\begin{definition}
Let $H$ and $G$ be subgroups of $\aut X^*$ and $n$ be a positive integer. We write $H^{2^n}\preceq G$ if $G$ contains $2^n$ copies of $H$ acting on the $2^n$ subtrees of level $n$.
\end{definition}
\begin{proposition}\label{prop_dimb}
Let $H$ and $G$ be subgroups of $\aut X^*$, such that $H^{2^n}\preceq G$. Then $\dim_HH\leq\dim_HG$.
\end{proposition}
\begin{proof}
It is straightforward that $H^{2^n}\preceq G\implies\abs{\quo{H}{(m-n)}}^{2^n}\leq\abs{\quo{G}{m}}$ for $m\geq n$. The conclusion follows.
\end{proof}

We now turn to the construction of a full-dimensional group.
Let $a_1=\sigma\in\aut X^*$ be the permutation exchanging the two maximal subtrees and let $a_n$ be defined recursively as
\begin{align*}
a_n=\pair<1,a_{n-1}>\sigma.
\end{align*}
It is straightforward that $a_n$ is of order $2^n$. It can be viewed as a finite-depth version of the familiar adding machine $t=\pair<1,t>\sigma$. The important thing is that $a_n$ acts as a full cycle on the $n$-th level vertices of the tree, and trivially below the $n$-th level.

Next, for any element $g\in\aut X^*$ and any word $w\in X^*$, we define $w*g\in\aut X^*$ as the element which acts as $g$ on the subtree $wX^*$, and trivially everywhere else. The following identity can be checked directly:
\begin{align*}
(w*g)^h=w^h*g^{(h@w)}.
\end{align*}
Let $g_1,\dots,g_n$ be any elements in $\aut X^*$. We define the element
\begin{align*}
\delta(g_1,\dots,g_n)=\prod_{i=0}^{n-1}(1^i0^{n-i})*g_{i+1}.
\end{align*}
Notice that $1^i0^{n-i}=(1^n)^{(a_n)^{2^i}}$, and that the product above can be taken in any order as the elements all commute (they act on different subtrees).
\begin{lemma}\label{lem_delta}
Let $G=\gen{g_1,\dots,g_n}$ and $H=\gen{a_n,\delta(g_1,\dots,g_n)}$. Then $(G')^{2^n}\preceq H'$.
\end{lemma}
\begin{proof}
The following equalities are immediate consequences of the definitions:
\begin{align*}
\delta(g_1,\dots,g_n)^{(a_n)^{k}}&=\prod_{i=0}^{n-1}(1^n)^{(a_n)^{2^i+k}}*g_{i+1},\\
\left[\delta(g_1,\dots,g_n),a_n^{k}\right]\left[\delta(g_1,\dots,g_n)^{-1},a_n^{k}\right]&=(1^n)^{(a_n)^{2^i}}*\left[g_{j+1},g_{i+1}\right]^{g_{j+1}^{-1}}.
\end{align*}
The second relation holds whenever there are integers $i,j\in\{0,\dots,n-1\}$ such that $2^i-2^j=k$.
These two equalities imply
\begin{align*}
w*[g_i,g_j]\in H'
\end{align*}
for all $w\in X^n$ and $i,j\in\{1,\dots,n\}$. This in turn implies
\begin{align*}
(G')^{2^n}\preceq H'.
\end{align*}
\end{proof}

Let $B_n=C_2^n$ be the direct product of $n$ copies of $C_2$. Let $\omega_1,\dots,\omega_n$ be a basis of $B_n^*$ and consider the spinal group $G_n$ defined through the sequence $\omega_1\dots\omega_n\omega_1\dots\omega_n\dots$. In other words $G_n=\gen{a,b_{(1,n)},\dots,b_{(n,n)}}$, with
\begin{align*}
a&=a_1=\sigma,\\
b_{(i,n)}&=\pair<1,b_{(i+1,n)}>&\text{for }i=1,\dots,n-1,\\
b_{(n,n)}&=\pair<a,b_{(1,n)}>.
\end{align*}
It follows from Corollary~\ref{cor_dim2} that $\dim_HG_n=1-\frac{3}{2^{n+1}}$. Define the elements
\begin{align*}
\tilde b_n=\delta(a,b_{(1,n-3)},\dots,b_{(n-3,n-3)},a_{n+1},\tilde b_{n+1})
\end{align*}
for $n\geq3$, and write $H_n=\gen{a_n,\tilde b_n}$.
\begin{theorem}
$H_n$ has Hausdorff dimension equal to $1$, for all $n\geq3$.
\end{theorem}
\begin{proof}
Lemma~\ref{lem_delta} yields
\begin{align*}
(H_{m+1}')^{2^m}\preceq H_m',&&(G_{m-3}')^{2^m}\preceq H_m'
\end{align*}
for all $m\geq3$. Thus $\dim_HH_m'\geq\dim_HG_{m-3}'$ by Proposition~\ref{prop_dimb}, and $\dim_HH_n'\geq\dim_HH_m'$ for all $3\leq n\leq m$. Proposition~\ref{prop_dima} allows us to state
\begin{align*}
\dim_HH_n\geq\dim_HH_n'\geq\dim_HG_m'=\dim_HG_m
\end{align*}
for all $3\leq n\leq m+3$. The last equality holds because $G_m$ is generated by $m$ elements of order $2$, so $\abs{G_m:G_m'}$ is finite. This yields $\dim_HH_n=1$ for all $n\geq3$, since $\dim_HG_m=1-\frac{3}{2^{m+1}}$.
\end{proof}

\begin{remark}
We could easily extend this construction by taking any sequence of finitely generated groups $G_n$ such that $\limsup(\dim_HG_n')=1$.
\end{remark}

\section{Proof of the main theorem}\label{sec_pr}

The remainder of the article is devoted to the proof of Theorem~\ref{th_main}.
Our first goal is to find a recursive formula for $\abs{\quo{G_\omega}{m}}$. We begin with a few simple but very useful lemmata.

Let $\pi:\aut\T\to(\quo{\aut\T}{1})\simeq C_2$ be the natural epimorphism.
\begin{lemma}[Folklore]
The map $\phi_n:\aut\T\to C_2$ given by $g\mapsto\prod_{w\in X^n}\pi(g@w)$ is an epimorphism for all $n\in\N$.
\end{lemma}
For a group $G\leq\aut\T$, we let $\stab[G]{n}$ be the subgroup of $G$ consisting of the elements that fix the first $n$ levels of the tree.
When $G=\aut\T$ we simply write $\stab{n}$.
For $v\in X^m$, we define $v\phi_n:\stab{m}\to C_2$ by $v\phi_n(g)=\prod_{w\in X^n}\pi(g@(vw))$.
\begin{corollary}
The map $v\phi_n:\stab{m}\to C_2$ is an epimorphism for all $v\in X^m$ and $m,n\in\N$.
\end{corollary}
\begin{proof}
This is straightforward because $\stab{m}\simeq(\aut\T)^{X^m}$.
\end{proof}

In the following we let $\gen{\quo{\omega}{m}}$ (resp. $\gen{\omega}$) designate the vector space spanned by $\omega_1,\dots,\omega_{s^{-1}(m)}$ (resp. $\omega_1,\omega_2,\dots$).

Consider $\psi\in\gen{\omega}$ and $x\in X$.
We define the homomorphisms $\overline\psi^x:\stab[G_\omega]{1}\to C_2$ as follows.
Write $\psi=\omega_{i_1}+\cdots+\omega_{i_k}$ where the $\omega_{i_j}$ are pairwise distinct and all appear at least once in $\omega$.
Let $(n_j+1)$ be the position of the first occurence of $\omega_{i_j}$ in $\omega$.
If $n_j>0$ for all $j$ then we define $\overline{\psi}^x=\sum_{j=1}^kx\phi_{n_j}$ (we implicitly identify $C_2$ with $\Z/2\Z$).
Otherwise if $n_1=0$ we set $\overline{\psi}^x=\bar x\phi_0 + \sum_{j=2}^kx\phi_{n_j}$, with $\bar x=1-x$.

By construction, $\overline{\psi}^x$ is a homomorphism.
What may be less obvious is the following lemma
(Notice that the group $\stab[G_\omega]{1}$ is generated by $\{b,b^a:b\in B\}$).
\begin{lemma}
Consider $\psi\in\gen{\omega}$.
Then the map $\overline{\psi}^0:\stab[G_\omega]{1}\to C_2$ (resp. $\overline{\psi}^1$) is the homomorphism induced by $b\mapsto\psi(b)$ and $b^a\mapsto1$ (resp. $b\mapsto1$ and $b^a\mapsto\psi(b)$) for $b\in B$.
In particular $\overline{\psi}^x$ is surjective.
\end{lemma}
\begin{proof}
This is easy to check in case $\psi$ appears in the sequence $\omega$.
The general case is just a linear combination of the terms of $\omega$.
\end{proof}

\begin{remark}
One can define $\overline{\psi}^x:(\quo{\stab[G_\omega]{1}}{m})\to C_2$ for any $\psi\in\gen{\quo{\omega}{m}}$ in the same way and the lemma still holds.
\end{remark}

To proceed further we need to define some specific subgroups of $G_\omega$. Let $\psi$ be an element of $B^*$. We define the subgroup $T_\omega(\psi)=\gen{\ker\psi}^{G_\omega}$, where the superscript designates normal closure in $G_\omega$. It should be noted that $T_\omega(\psi)\leq\stab[G_\omega]{1}$ for every $\psi\in B^*$. We will now state and prove two technical lemmata, which lead to Proposition~\ref{prop1}.

\begin{lemma}\label{tphi}
Let $\omega=\omega_1^{a_1}\omega_2^{a_2}\dots$ be in syllable form; let $m\geq1$ be an integer; let $\psi_0$ be a non-trivial element in $B^*$ such that $\psi_0\neq\omega_1$. Then
\begin{align*}
  \log_2\abs{(\quo{G_\omega}{m})/(\quo{T_\omega(\psi_0)}{m})}=
  \begin{cases}
    1 &\text{if }m\leq s_k+1,\\
    3 &\text{if }m>s_k+1,
  \end{cases}
\end{align*}
where $k$ is the greatest integer such that $\psi_0$ is linearly independent from $\omega_1,\dots,\omega_k$.
\end{lemma}
\begin{proof}
Case $m\leq s_k+1$.
Let $\{\psi_1,\dots,\psi_\lambda\}$ be a basis of $\gen{\quo{\omega}{m}}$.
Since $\psi_0\not\in\gen{\quo{\omega}{m}}$, the set $\{\psi_0,\dots,\psi_\lambda\}$ is a basis of some subspace of $B^*$.
Let $\{b_0,\dots,b_\lambda\}\subset B$ be a dual basis, i.e. $b_0,\dots,b_\lambda$ satisfy $\psi_i(b_j)=\delta_{ij}$ for all $i,j\in\{0,\dots,\lambda\}$.

We have $b_1,\dots,b_\lambda\in\ker\psi_0$, so $b_1,\dots,b_\lambda\in T_\omega(\psi_0)$.
Since $(\quo{\stab[G_\omega]{1}}{m})$ is generated as a normal subgroup by the images of $b_1,\dots,b_\lambda$, we have $(\quo{T_\omega(\psi_0)}{m})=(\quo{\stab[G_\omega]{1}}{m})$.
Therefore $(\quo{G_\omega}{m})/(\quo{T_\omega(\psi_0)}{m})=C_2$.

Case $m>s_k+1$.
Let $\{\psi_0,\dots,\psi_\lambda\}$ be a basis of $\gen{\omega}$ and let $\{b_0,\dots,b_\lambda\}$ be a dual basis.
Write $H=\stab[G_\omega]{1}=\gen{b_0,\dots,b_\lambda}^{G_\omega}$.
Then obviously $(\quo{G_\omega}{m})/(\quo{H}{m})=C_2$.

We will now prove that $(\quo{H}{m})/(\quo{T_\omega(\psi_0)}{m})=C_2\times C_2$.
This group is generated by the images of $b_0$ and $b_0^a$.
A straightforward computation shows that $b_0^2=[b_0,b_0^a]=1$.
Therefore $(\quo{H}{m})/(\quo{T_\omega(\psi_0)}{m})$ is a quotient of $C_2\times C_2$.
Consider the map $\Psi:g\mapsto\left(\overline{\psi_0}^0(g),\overline{\psi_0}^1(g)\right)$.
It is a surjective group homomorphism $(\quo{H}{m})\to C_2\times C_2$, but $\Psi(T_\omega(\psi_0))$ is trivial.
Therefore $\Psi(\quo{T_\omega(\psi_0)}{m})$ has index $4$ in $\Psi(\quo{H}{m})$, which finishes the proof.
\end{proof}

\begin{lemma}[Case $\phi=\omega_1$]\label{t1}
Let $\omega=\omega_1^{a_1}\omega_2^{a_2}\dots$ be in syllable form and set $m\geq1$. Then
\begin{align*}
  \log_2\abs{(\quo{G_\omega}{m})/(\quo{T_\omega(\omega_1)}{m})}=
  \begin{cases}
    1 &\text{if }m=1,\\
    m+1 &\text{if }1<m\leq a_1+1,\\
    a_1+2 &\text{if }a_1+1<m\leq s_k+1,\\
    a_1+3 &\text{if }s_k+1<m,
  \end{cases}
\end{align*}
where $k$ is the greatest integer such that $\omega_1$ is linearly independent from $\omega_2,\dots,\omega_k$.
\end{lemma}
\begin{proof}
The case $m=1$ is very simple because $(\quo{T_\omega(\omega_1)}{1})$ is the trivial group and $(\quo{G_\omega}{1})=C_2$.

Define $H=\stab[G_\omega]{a_1+1}$.
If $1<m\leq a_1+1$, since $T_\omega(\omega_1)\leq H$, we know that $(\quo{T_\omega(\omega_1)}{m})$ is trivial.
It is clear that $(\quo{G_\omega}{m})$ is isomorphic to a dihedral group of order $2^{m+1}$ because $(\quo{G_\omega}{m})$ is generated by two involutions $a$ and $b_1$, and $ab_1$ has order $2^m$ in $(\quo{G_\omega}{m})$.

Case $a_1+1<m\leq s_k+1$.
Let $\{\omega_{i_1},\dots,\omega_{i_\lambda}\}$ be a basis of $\gen{\omega}$, with $\omega_{i_1}=\omega_1$.
Let $\{b_1,\dots,b_\lambda\}$ be a dual basis.
It is readily checked that the group $\gen{a,b_1}$ is dihedral of order $2^{a_1+3}$, and a straightforward computation shows that $a^2=b_1^2=1$ and $(ab_1)^{2^{a_1+1}}\in\stab[G_\omega]{s_k+1}$.
Since $(\quo{G_\omega}{m})/(\quo{T_\omega(\omega_1)}{m})$ is generated by the images of $a$ and $b_1$, we conclude that this group is dihedral of order $2^{a_1+2}$.

Finally, if $s_k+1<m$, it is sufficient to prove $(\quo{H}{m})/(\quo{T_\omega(\omega_1)}{m})\simeq C_2$.
By the above we know that this group is generated by the image of $(ab_1)^{2^{a_1+1}}$, and that this element is of order $2$.
Hence $(\quo{H}{m})/(\quo{T_\omega(\omega_1)}{m})$ is a quotient of $C_2$.
Express $\omega_1$ as a linear combination of $\omega_2,\dots,\omega_{k+1}$: $\omega_1=\omega_{i_1}+\cdots+\omega_{i_\lambda}$.
Let $n_j$ be the position of the last occurence of $\omega_{i_j}$ in $(\quo{\omega}{m})$.
Fix $v=0^{a_1+1}$ and define $\overline{\omega_1}:(\quo{H}{m})\to C_2$ by $\overline{\omega_1}=\sum_{j=1}^\lambda v\phi_{n_j}$.
Then $\overline{\omega_1}$ is surjective but $\overline{\omega_1}(\quo{T_\omega(\omega_1)}{m})$ is trivial.
This completes the proof.
\end{proof}

\begin{proposition}\label{prop1}
Let $\omega=\omega_1^{a_1}\omega_2^{a_2}\dots$ be in syllable form and set $m>a_1$. Then
\begin{align*}
  \log_2\abs{\quo{G_{\omega}}{m}}=
    2+a_1+\delta(m)+2^{a_1}\bigl(\log_2\abs{\quo{G_{\sigma^{a_1}\omega}}{(m-a_1)}}-2\delta(m)-1\bigr),
\end{align*}
with
\begin{align*}
  \delta(m)=\begin{cases}0&\text{if }\omega_1\text{ is linearly independent from }\omega_{2},\dots,\omega_{s^{-1}(m)},\\1&\text{otherwise.}\end{cases}
\end{align*}
\end{proposition}
\begin{proof}
For $m>a_1$, Lemma~\ref{t1} gives
\begin{align*}
  \log_2\abs{(\quo{G_\omega}{m})/(\quo{T_\omega(\omega_1)}{m})}=
  \begin{cases}
    a_1+2 &\text{if }m\leq s_k+1,\\
    a_1+3 &\text{if }m>s_k+1,
  \end{cases}
\end{align*}
where $k$ is the greatest integer such that $\omega_1$ is linearly independent from $\omega_2,\dots,\omega_k$. We can rewrite this equation as
\begin{align}\label{eqgt1}
  \log_2\abs{\quo{G_\omega}{m}}=
    a_1+2+\delta(m)+\log_2\abs{\quo{T_\omega(\omega_1)}{m}}.
\end{align}
Next, iteration of the relation $T_\omega(\omega_1)=T_{\sigma\omega}(\omega_1)\times T_{\sigma\omega}(\omega_1)$ gives
\begin{align*}
  T_\omega(\omega_1)=\underbrace{T_{\sigma^{a_1}\omega}(\omega_1)\times\cdots\times T_{\sigma^{a_1}\omega}(\omega_1)}_{2^{a_1}},
\end{align*}
therefore
\begin{align}\label{eqgt2}
  \log_2\abs{\quo{T_\omega(\omega_1)}{m}}=2^{a_1}\log_2\abs{\quo{T_{\sigma^{a_1}\omega}(\omega_1)}{(m-a_1)}}.
\end{align}
But Lemma~\ref{tphi} yields
\begin{align*}
  \log_2\abs{(\quo{G_{\sigma^{a_1}\omega}}{(m-a_1)})/(\quo{T_{\sigma^{a_1}\omega}(\omega_1)}{(m-a_1)})}=
  \begin{cases}
    1 &\text{if }m\leq s_{k'}+1,\\
    3 &\text{if }m>s_{k'}+1,
  \end{cases}
\end{align*}
where $k'$ is the greatest integer such that $\omega_1$ is linearly independent from $\omega_2,\dots,\omega_{k'}$, i.e. $k=k'$. Therefore we can rewrite the preceding equation as
\begin{align}\label{eqgt3}
  \log_2\abs{\quo{G_{\sigma^{a_1}\omega}}{(m-a_1)}}=
    1+2\delta(m)+\log_2\abs{\quo{T_{\sigma^{a_1}\omega}(\omega_1)}{(m-a_1)}}.
\end{align}
Equations~\eqref{eqgt1}, \eqref{eqgt2} and~\eqref{eqgt3} give the result.
\end{proof}

Now the technical part is over, and the following statements and their proof, including the proof of Theorem~\ref{th_main}, are easy consequences of what has been shown above.
\begin{proposition}\label{prop_6.8}
Let $\omega=\omega_1^{a_1}\omega_2^{a_2}\dots$ be in syllable form and consider $\lambda\in\N$ and $m>s_\lambda$. Then
\begin{multline*}
  \log_2\abs{\quo{G_{\omega}}{m}}=
3+2^{s_0}(1+a_1+\delta_1-2\delta_0)+\cdots+2^{s_{\lambda-1}} (1+a_\lambda+\delta_\lambda-2\delta_{\lambda-1}) \\+ 2^{s_\lambda}(\log_2\abs{\quo{G_{\sigma^{s_\lambda}\omega}}{(m-s_\lambda)}}-2\delta_\lambda-1),
\end{multline*}
with
\begin{align*}
  \delta_j=\delta_j(m)=\begin{cases}0&\text{if }\omega_j\text{ is linearly independent from }\omega_{j+1},\dots,\omega_{s^{-1}(m)},\\1&\text{otherwise.}\end{cases}
\end{align*}
(We set $\delta_0=1$ and $s_0=0$).
\end{proposition}
\begin{proof}
This follows directly from $\lambda$ applications of Proposition~\ref{prop1}.
\end{proof}

\begin{corollary}\label{cor_6.9}
Let $\lambda$ be such that the space spanned by $\omega_{\lambda+1},\dots,\omega_{s^{-1}(m)}$ contains all $\omega_j$ with $1\leq j\leq s^{-1}(m)$. Then
\begin{align*}
  \log_2\abs{\quo{G_{\omega}}{m}}=
3+\Sigma_\lambda+2^{s_\lambda}(\log_2\abs{\quo{G_{\sigma^{s_\lambda}\omega}}{(m-s_\lambda)}}-3),
\end{align*}
with
\begin{align*}
  \Sigma_\lambda=2^{s_0}a_1+\cdots+2^{s_{\lambda-1}} a_\lambda.
\end{align*}
\end{corollary}

\begin{lemma}\label{lem_dim_fall}
Let $\omega=\omega_1^{a_1}\omega_2^{a_2}\dots$ be in syllable form. Given $m>a_1+1$, if $\lambda$ is the smallest integer such that $\omega_\lambda$ is linearly independent from $\omega_{\lambda+1},\dots,\omega_{s^{-1}(m)}$, then
\begin{align*}
  \log_2\abs{\quo{G_{\omega}}{m}}=
3+\Sigma_\lambda+2^{s_\lambda+1}-2^{s_{\lambda-1}}+2^{s_\lambda}(\log_2\abs{\quo{G_{\sigma^{s_\lambda}\omega}}{(m-s_\lambda)}}-3).
\end{align*}
\end{lemma}
\begin{proof}
This is just a consequence of Proposition~\ref{prop_6.8} and Corollary~\ref{cor_6.9}.
\end{proof}
\begin{remark}\label{rem1}
If $1<m\leq a_1+1$ then $(\quo{G_{\omega}}{m})$ is just a dihedral group of order $m+1$, whence
\begin{align*}
  \log_2\abs{\quo{G_{\omega}}{m}}=m+1.
\end{align*}
\end{remark}

We are naturally led to the following proposition.
\begin{proposition}\label{prop:sizegmodm}
Consider $\omega=\omega_1^{a_1}\omega_2^{a_2}\dots$ in syllable form with $\dim(\quo{\omega}{m})=n$. Then \begin{align}\label{quofin}
  \log_2\abs{\quo{G_{\omega}}{m}}=
3+\Sigma_{s^{-1}(m)-1}+\sum_{i=2}^{n-1}(2^{s_{\lambda_i}+1}-2^{s_{\lambda_i-1}})-2^{s_{\lambda_1-1}}
+2^{s_{\lambda_1}}(m-s_{\lambda_1}),
\end{align}
where $\lambda_j$ is the smallest integer such that $\dim(\quo{\sigma^{s_{\lambda_j}}\omega}{(m-s_{\lambda_j})})=j$, for each $j\in\{1,\dots,n-1\}$.
\end{proposition}
\begin{proof}
We simply apply the previous lemma $n-1$ times to obtain
\begin{align*}
  \log_2\abs{\quo{G_{\omega}}{m}}=
3+\Sigma_{\lambda_1}+\sum_{i=1}^{n-1}(2^{s_{\lambda_i}+1}-2^{s_{\lambda_i-1}})
+2^{s_{\lambda_1}}(\log_2\abs{\quo{G_{\sigma^{s_{\lambda_1}}\omega}}{(m-s_{\lambda_1}})}-3).
\end{align*}
Next, we have $\dim(\quo{\sigma^{s_{\lambda_1}}\omega}{(m-s_{\lambda_1})})=1$ and $\lambda_1=s^{-1}(m)-1$.
Remark~\ref{rem1} yields
\begin{align*}
  \log_2\abs{\quo{G_{\omega}}{m}}=
3+\Sigma_{s^{-1}(m)-1}+\sum_{i=1}^{n-1}(2^{s_{\lambda_i}+1}-2^{s_{\lambda_i-1}})
+2^{s_{\lambda_1}}(m-s_{\lambda_1}-2).
\end{align*}
The result is obtained by extracting the first term of the sum.
\end{proof}

We are now ready to prove Theorem~\ref{th_main}, which we restate here.

\begin{theorem}
Consider $\omega=\omega_1^{a_1}\omega_2^{a_2}\dots$ in syllable form with $\diminf(\omega)=n\geq2$. The Hausdorff dimension of $G_\omega$ is computed as
\begin{align*}
  \dim_HG_{\omega}=\frac{1}{2}\liminf_{k\to\infty}\left(\frac{\Sigma_k}{2^{s_k}}
+\frac{1}{2^{s_k}}\sum_{i=2}^{n-1}2^{s_{\lambda_i}}\left(2-\frac{1}{2^{a_{\lambda_i}}}\right)
+\frac{1}{2^{a_k}}\left(1-\frac{1}{2^{a_{k-1}}}\right)\right),
\end{align*}
where for each $i\in\{2,\dots,n-1\}$ we let $\lambda_i(k)$ be the smallest integer such that
\begin{align*}
  \dim(\quo{\sigma^{s_{\lambda_i(k)}}\omega}{(s_{k+1}-s_{\lambda_i(k)})})=i.
\end{align*}
\end{theorem}
\begin{proof}
Starting with Equation~\eqref{quofin}, we write $k=s^{-1}(m)$. Recalling $\lambda_1=k-1$, we compute
\begin{multline*}
  \frac{\log_2\abs{\quo{G_{\omega}}{m}}}{2^m}=\\\frac{1}{2^{m-s_{k}}}\left(
\frac{3}{2^{s_{k}}}+\frac{\Sigma_{k-1}}{2^{s_{k}}}
+\frac{1}{2^{s_{k}}}\sum_{i=2}^{n-1}(2^{s_{\lambda_i}+1}-2^{s_{\lambda_i-1}})
-\frac{1}{2^{s_{k}-s_{k-2}}}+\frac{1}{2^{a_{k}}}(m-s_{k-1})\right).
\end{multline*}
If we fix $k$ and consider $m$ such that $s_{k-1}+1<m\leq s_k+1$, the $\lambda_i$'s do not depend on $m$. We easily check that the expression is minimal when $m=s_{k}+1$. Therefore
\begin{align*}
  \dim_HG_{\omega}=\frac{1}{2}\liminf_{k\to\infty}\left(\frac{\Sigma_k}{2^{s_k}}
+\frac{1}{2^{s_k}}\sum_{i=2}^{n-1}(2^{s_{\lambda_i}+1}-2^{s_{\lambda_i-1}})
-\frac{1}{2^{a_k+a_{k-1}}}
+\frac{1}{2^{a_k}}\right).
\end{align*}
The result follows.
\end{proof}

\bibliography{../bibliography/trees}
\end{document}